\documentclass[12pt]{article}
\usepackage{}
\usepackage{mathrsfs}
\usepackage{amsfonts}

\usepackage{amsthm}
\usepackage{indentfirst,latexsym,mathrsfs,amsmath,amssymb,fancyhdr,setspace}
\usepackage{graphicx}
\usepackage{color}

\numberwithin{equation}{section}

\newcommand{\be}{\begin{equation}\label}
\newcommand{\ee}{\end{equation}}

\newcommand{\beaa}{\begin{eqnarray}}
\newcommand{\bea}{\begin{eqnarray}\label}
\newcommand{\eea}{\end{eqnarray}}

\newcommand{\bas}{\begin{eqnarray*}}
\newcommand{\eas}{\end{eqnarray*}}

\newcommand{\Om}{\Omega}

\newcommand{\di}{\textrm{div}}

\newcommand{\cd}{\cdot}
\newcommand{\nn}{\nonumber}

\newcommand{\ep}{\epsilon}
\newcommand{\epz}{\varepsilon}

\newcommand{\f}{\frac}
\newcommand{\p}{\partial}

\newcommand{\na}{\nabla}
\newcommand{\del}{\delta}
\newcommand{\Del}{\Delta}

\newcommand{\bu}{\mathbf u}
\newcommand{\bv}{\mathbf v}
\newcommand{\bV}{\mathbf V}
\newcommand{\bw}{\mathbf w}
\newcommand{\bh}{\mathbf h}

\newcommand{\R}{\mathbb{R}}

\newtheorem{thm}{Theorem}[section]
\newtheorem{lem}{Lemma}[section]
\newtheorem{coro}{Corollary}[section]
\newtheorem{dnt}{Definition}[section]
\newtheorem{remark}{Remark}[section]

\newtheorem{prop}{Proposition}[section]

\title{Relative Entropy Method for the relaxation limit of Hydrodynamic models}
\date{}

\author{\small{JOS\'{E} A. CARRILLO, YINGPING PENG AND ANETA WR\'{O}BLEWSKA-KAMI\'{N}SKA}}

\begin{document}
\maketitle
{\small
{\textbf{Abstract:}
We show how to obtain general nonlinear aggregation-diffusion models, including Keller-Segel type models with nonlinear diffusions, as relaxations from nonlocal compressible Euler-type hydrodynamic systems via the relative entropy method. We discuss the assumptions on the confinement and interaction potentials depending on the relative energy of the free energy functional allowing for this relaxation limit to hold. We deal with weak solutions for the nonlocal compressible Euler-type systems and strong solutions for the limiting aggregation-diffusion equations. Finally, we show the existence of weak solutions to the nonlocal compressible Euler-type systems satisfying the needed properties for completeness sake.
\\

{\textbf{Keywords:} relative entropy method, aggregation-diffusions, nonlocal hydrodynamics, relaxation limit.


\section{Introduction}

In this work, we consider the following compressible Euler-type systems of equations of the form
\begin{equation}\label{Eqs-1}
\begin{split}
\p_t\rho   +    \di_x (\rho\bu)   &=0, \\
\p_t(\rho\bu)  +  \di_x(\rho\bu\otimes\bu)  &=   -\f{1}{\epz}\rho\na_{x}\f{\del\mathcal{E}(\rho)}{\del\rho}
     -  \f{1}{\epz}\rho\bu
\end{split}
\end{equation}
in the time-spatial domain $(0,T)\times\Om$, where $\rho(t):\Om\rightarrow\R_{+}$ for $t\geq0$ is the density obeying the equation of conservation of mass,
$\bu(t):\Om\rightarrow\R^d$ for $t\geq0$ is the velocity of fluid and the product $\rho\bu$ denotes the momentum flux. Here the functional $\mathcal{E}(\rho): L^1_+(\R^d)\rightarrow\R$ is the free energy functional defined on mass densities by
\begin{align}\label{eq:freeenergy}
\mathcal{E}(\rho)  = \int_{\mathbb{R}^d} h(\rho)dx +\int_{\mathbb{R}^d}  \Phi( x ) \rho\, d x + \frac{C_k}{2} \int_{\mathbb{R}^d}   (K\ast\rho) \rho \,d x,
\end{align}
with $h(\rho)$ describing the entropy part or internal energy of the system, and $\f{\del\mathcal{E}(\rho)}{\del\rho}$ stands for its variational derivative, given by
\begin{align}\label{E-h-1}
\f{\del\mathcal{E}(\rho)}{\del\rho} = h'(\rho)+\Phi+C_k (K\ast\rho)\,.
\end{align}
Here, $C_k$ is a positive constant measuring the strength of the interaction, $K(x):\R^d\rightarrow\R$ is the interaction potential depicting the nonlocal forces which usually manifest as repulsion or attraction between particles, which is assumed to be symmetric, and
$\Phi(x):\Om\rightarrow\R$ is a confinement potential. We refer to \cite{CMV2,CMV1,V} for a general introduction to these free energies, to \cite{VC} for their applications in Keller-Segel type models, and more general models in Density Functional Theory as discussed in \cite{GPK}.
Finally, the term $-\f{1}{\epz}\rho\bu$ on the left-hand-side of \eqref{Eqs-1} is responsible for a damping force with frictional coefficient $\f{1}{\epz}$ in order to look at the so-called overdamped limit.

In this work, we consider $\Om\subset\R^d$ to be any smooth, connected, open set. The no-flux boundary condition for $\bu$ (i.e. $\bu\cd\nu=0$, $\nu$ denotes an outer normal vector to $\partial \Omega$)) or periodic boundary condition are assumed if $\Om$ is a bounded domain or $\Om=\mathbb{T}^d$ is periodic domain. We also extend $\rho$ by zero when $\Om$ is bounded in order that we are able to define properly $K\ast\rho$ on $\R^d$. The main objective of this work is to deduce the following equilibrium equation
\begin{align}\label{equi-eq-1}
\p_t\bar{\rho}   =  \di_x\left(\bar{\rho}\na_{x}\f{\del\mathcal{E}(\bar{\rho})}{\del\rho} \right)
\end{align}
by taking the overdamped limit $\epz\rightarrow0$ in system \eqref{Eqs-1} under the framework of relative entropy method. This method is an efficient mathematical tool for establishing the limiting processes and stabilities among thermomechanical theories, see \cite{JAC-1,CCT,CMD-1,CMD-2,RDJ,JG,CL-2,CL-1} for instance. With the various choices of the functional $\mathcal{E}(\rho)$, the corresponding models spanned from the system of isentropic gas dynamics and variants of the Euler-Poisson system \cite{QJ,TL,YP} leading to the porous medium equation and nonlinear aggregation-diffusion equations in the overdamped limit, see \cite{CG,H1,H2,H3,JR,MM} and references therein. More general forms of free energies with higher order terms in derivatives have also been used in the literature leading to the equations of quantum hydrodynamics \cite{PA-2,PA-1}, the models for phase transitions \cite{DB,XR}, and the dispersive Euler-Korteweg equations \cite{JED}.

In this work, we only consider the functional $\mathcal{E}(\rho)$ defined by \eqref{eq:freeenergy} with variation given by \eqref{E-h-1}
where $h(\rho)$ and a pressure function denoted by $p(\rho)$ are linked by the thermodynamic consistency relations
\begin{align}\label{h-p-1}
\rho h''(\rho)=p'(\rho),\quad \rho h'(\rho)=p(\rho)  +  h(\rho).
\end{align}
In this case, we observe that \eqref{Eqs-1} reduces to
\begin{equation}\label{Eqs}
\begin{split}
\p_t\rho   +    \di_x (\rho\bu)   &=0  , \\
\p_t(\rho\bu)  +  \di_x(\rho\bu\otimes\bu) +  \f{1}{\epz}\na_x p(\rho)
&=  -  \f{C_k}{\epz}(\na_x K\ast\rho)\rho  -  \f{1}{\epz}\rho\bu
         -  \f{1}{\epz}\rho\na_x\Phi
\end{split}
\end{equation}
and  \eqref{equi-eq-1} is equivalent to
\begin{align}\label{equi-eq}
\p_t\bar{\rho}   =  \Del_x p(\bar{\rho})   +   C_k\di_x((\na_x K\ast\bar{\rho})\bar{\rho})  +   \di_x(\bar{\rho}\na_x\Phi),
\end{align}
consequently, our goal concerning the relaxation limit from \eqref{Eqs-1} to \eqref{equi-eq-1} is equivalent to considering the relaxation limit from \eqref{Eqs} to \eqref{equi-eq}. In particular, for the power-law pressure $p(\rho)=\rho^m$, the internal energy $h(\rho)$ takes the form
\begin{equation*}
h(\rho)=
\left\{
\begin{split}
&\f{1}{m-1}\rho^{m},\quad m>1 ,  \\
&\rho \log\rho,   \quad m=1.
\end{split}
\right.
\end{equation*}
We will deal with slightly more general internal energy functions. For this reason, we introduce the notation
\begin{equation}\label{h-1}
h_m(\rho)=
\left\{
\begin{split}
&k_1\rho\log\rho ,\quad m=1,   \\
&\f{k_2}{m-1}\rho^{m},\quad 1<m\le2,   \\
&\f{k_3}{m-1}\rho^{m}   + o(\rho^m) \quad \textrm{as}\,\, \rho\rightarrow+\infty, \quad m>2, \ h_m \in C[0,+\infty)\cap C^2(0,+\infty), \ h_m^{''}(\varrho)>0 %
\end{split}
\right.
\end{equation}
for some positive constants $k_1$, $k_2$ and $k_3$. For $m>2$, we assume that the function $o(\rho^m)$ is chosen to satisfy that for some constant $A>0$,
\begin{equation}\label{pA}
|p''(\rho)|  \le  A\f{p'(\rho)}{\rho}  \qquad \forall\rho>0,
\end{equation}
where $p(\rho)$ is determined by $h_m(\rho)$ via \eqref{h-p-1}. For simplicity, we will drop the dependence on $m$ of $h(\rho)$ in the sequel.

We can formally obtain that  weak solutions $(\rho,\rho\bu)$ of the system \eqref{Eqs} satisfy a standard weak form of total energy dissipation. Indeed, multiplying \eqref{Eqs}$_2$ with $\bu$, using \eqref{Eqs}$_1$ and integrating the resulting equation over $\Om$, provided no-flux boundary condition for $\bu$ (i.e. $\bu\cd\nu=0$) is valid when $\Om\subset\R^d$ is a bounded domain, one derives
 \begin{align}\label{energy-1}
 \f{d}{dt}\int_{\Om}\left(\f{1}{\epz}h(\rho)  +  \f12\rho|\bu|^2 +   \f{C_k}{2\epz}(K\ast\rho)\rho
    +   \f{1}{\epz}\rho\Phi \right)dx +    \f{1}{\epz}\int_{\Om}\rho|\bu|^2dx      =0
 \end{align}
in the sense of distributions, where we have used the first relation in \eqref{h-p-1}.

In order to obtain the free energy dissipation for \eqref{equi-eq} and further to compare its strong solution with the  weak solution of \eqref{Eqs}, we define
\begin{align}\label{bar-m}
\bar{m}=\bar{\rho}\bar{\bu}= - \na_x p(\bar{\rho})
   -  C_k(\na_x K\ast\bar{\rho})\bar{\rho}   -   \bar{\rho}\na_x\Phi
\end{align}
and rewrite \eqref{equi-eq} as
\begin{equation}\label{Eqs-bar-1}
\begin{split}
\p_t\bar{\rho}   +   \di_{x}(\bar{\rho}\bar{\bu})  &=0,     \\
\p_t(\bar{\rho}\bar{\bu})   +   \di_{x}(\bar{\rho}\bar{\bu}\otimes\bar{\bu})
+\f{1}{\epz}\na_x p(\bar{\rho})
&=     -  \f{C_k}{\epz}(\na_x K\ast\bar{\rho})\bar{\rho}  -   \f{1}{\epz}\bar{\rho}\bar{\bu}   -\f{1}{\epz}\bar{\rho}\na_x\Phi   +   \bar{e},
\end{split}
\end{equation}
where $\bar{e}:=\p_t(\bar{\rho}\bar{\bu})  +   \di_x(\bar{\rho}\bar{\bu}\otimes\bar{\bu})$. In a similar way as for \eqref{energy-1}, we obtain the free energy dissipation for $(\bar{\rho},\bar{\rho}\bar{\bu})$ in the following form
\begin{align}\label{energy-2}
\f{d}{dt}\int_{\Om}\left(\f{1}{\epz}h(\bar{\rho})   +   \f{C_k}{2\epz}(K\ast\bar{\rho})\bar{\rho}
    +   \f{1}{\epz}\bar{\rho}\Phi\right) dx +    \f{1}{\epz}\int_{\Om}\bar{\rho}|\bar{\bu}|^2dx
    =  0,
\end{align}
where we have also assumed that no-flux boundary condition for $\bar{\bu}$ (i.e. $\bar{\bu}\cd\nu=0$) holds, when $\Om\subset\R^d$ is a bounded domain. Notice that this is the well-known dissipation property for gradient flows of the form \eqref{equi-eq-1}, see \cite{CMV2,CMV1,V} for instance.

For notational simplicity, we define the relative quantity $h(\rho|\bar{\rho})$ here by the difference between $h(\rho)$ and the linear part of the Taylor expansion around $\bar{\rho}$ as $h(\rho|\bar{\rho}):=h(\rho)  -  h(\bar{\rho})  -  h'(\bar{\rho})(\rho-\bar{\rho})$, and denote
\begin{align}\label{varphi-1}
\Theta(t):=
\f{1}{\epz}\int_{\Om}h(\rho|\bar{\rho})dx  +   \f12\int_{\Om}\rho|\bu-\bar{\bu}|^2dx + \f{C_k}{2\epz}\int_{\Om}(\rho-\bar{\rho})(K\ast(\rho-\bar{\rho}))dx,
\end{align}
which potentially measures the distance between the two solutions $(\rho,\rho\bu)$ and $(\bar{\rho},\bar{\rho}\bar{\bu})$. Indeed, assuming that the exponent of the pressure function satisfies
\begin{align}\label{restrict-m}
m\ge 2-\f{2}{d},\textrm{ for} \,\, d\ge2,
\end{align}
then the function $\Theta(t)$ provides a measure to the distance between $(\rho,\rho\bu)$ and $(\bar{\rho},\bar{\rho}\bar{\bu})$ in the relaxation limit as we will show below. The restrictions in \eqref{restrict-m} are due to the use of Hardy-Littlewood-Sobolev-type (HLS) inequalities. HLS inequalities are also essential for establishing the existence of global-in-time weak solutions to Keller-Segel systems for general initial data, see \cite{AB,VC,JAC,JW,YS} and references therein.
\begin{remark}\label{remark-2}
We should always keep in mind that whenever we deal with the equality case in \eqref{restrict-m}, the mass of our system \eqref{equi-eq} should be suitably smaller than a threshold value, called the critical mass, in order to deal without finite time blow-up problems, otherwise we can assume that time is small enough and deal with local in time solutions before the blow-up happens. For strict inequalities, we do not have any restrictions on the mass.
\end{remark}

We now recall the definition of weak solutions to \eqref{Eqs} we deal with in this work.

\begin{dnt}\label{def2}
$(\rho,\rho\bu)$ with $\rho\in C([0,T); L^1(\Om)\cap L^m(\Om))$, $\rho\geq0$ and $\rho|\bu|^2\in L^{\infty}(0,T;L^1(\Om))$ is a weak solution of \eqref{Eqs} if
\begin{itemize}
  \item $(\rho,\rho\bu)$ satisfies the weak form of \eqref{Eqs};
  \item $(\rho,\rho\bu)$ satisfies \eqref{energy-1} in the sense of distributions:
  \begin{align}\label{weak-en-1}
  &-\int_0^{\infty}\int_{\Om}\big(\f{1}{\epz}h(\rho) + \f12\rho|\bu|^2 + \f{C_k}{2\epz}(K\ast\rho)\rho + \f{1}{\epz}\rho\Phi\big)\dot{\theta}(t)dxdt
    +   \f{1}{\epz}\int_0^{\infty}\int_{\Om}\rho|\bu|^2\theta(t)dxdt   \nn \\
  &\quad =  \int_{\Om}\big(\f{1}{\epz}h(\rho) + \f12\rho|\bu|^2 + \f{C_k}{2\epz}(K\ast\rho)\rho + \f{1}{\epz}\rho\Phi\big)\big|_{t=0}\theta(0)dx
  \end{align}
  for any non-negative $\theta\in W^{1,\infty}[0,\infty)$ compactly supported on $[0,\infty)$;
  \item $(\rho,\rho\bu)$ satisfies the properties:
  \begin{align*}
         \int_{\Om}\rho(t,x)\, dx=M<\infty, \mbox{ for a.e. } t>0.
  \end{align*}
  \begin{align*}
  \sup_{t\in(0,T)}\int_{\Om}\left(\f{1}{\epz}h(\rho) + \f12\rho|\bu|^2 + \f{C_k}{2\epz}(K\ast\rho)\rho + \f{1}{\epz}\rho\Phi\right)dx<\infty.
  \end{align*}
\end{itemize}
\end{dnt}

Our main result is stated as follows.
\begin{thm}\label{th-1}
Let $T>0$ and $m\ge1$ be fixed. Let the confinement potential $\Phi(x)$ be bounded from below in $\Om$ and $p(\rho)$ be defined through \eqref{h-p-1} and \eqref{h-1} and let the interaction potential be symmetric. Suppose that $C_k$ is suitably small and $(\rho,\rho\bu)$ is a weak solution of \eqref{Eqs} in the sense of Definition \ref{def2} with $\rho>0$, and $(\bar{\rho},\bar{\rho}\bar{\bu})$ is a smooth solution of \eqref{equi-eq} with $\bar{\rho}>0$, $\bar{\bu} \in L^\infty(0,T;W^{1,\infty}(\Omega))\cap L^\infty(0,T; L^{\frac{m}{2(m-1)}}(\Om))$,  and $\bar{e}$ bounded.
Let $\Om$ be any smooth, connected, open subset in $\R^d$. Assume one of the following conditions hold:\\
\\
(i) $2-\f{2}{d}\le m \le 2$ with $d\ge2$
and the interaction potential $K$ satisfies
$K\in L^{\f{m}{2(m-1)}}(\Om)\cap W^{1,\infty}(\Om)$,  \\
  \\
(ii) $\Om=\mathbb{T}^d$ or $\Om$ is a bounded domain in $\R^d$, $m\ge 2-\f{2}{d}$ with $d\ge2$, $\bar{\rho}\in I=[\underline{\del},\overline{\del}]$ with $\underline{\del}>0$ and $\overline{\del}<\infty$ and the interaction potential $K$ satisfies
$ K\in L^p(\Om)\cap W^{1,\infty}(\Om)\,\, (1<p<\infty)$.   \\
  \\
Then the following stability estimate
\begin{align*}
\Theta(t)  \le C(\Theta(0)   +   \epz),  \quad t\in[0,T]
\end{align*}
holds, where $C$ is a positive constant depending only on $T$, possibly $I$, $\bar{\rho}$ and its derivatives. Moreover, if $\Theta(0)\rightarrow0$ as $\epz\rightarrow0$, then
\begin{align*}
\lim_{\epz\to 0} \sup_{t\in[0,T]}\Theta(t) = 0.
\end{align*}
\end{thm}

Notice that we may need more regular assumptions on the interaction potential $K$ and the confinement potential $\Phi$ in order to prove the existence of solutions to our systems. We will point out, in Section 3, the specific restrictions on $K$ and $\Phi$ when we show the existence of weak solutions to the system \eqref{Eqs} on two or three dimensional bounded domains. Otherwise, we just assume that $K$ and $\Phi$ are as regular as we need.

The outline of this paper is as follows. In Section 2, we first review how to obtain the relative entropy inequality for our system using the notion of weak solution in Definition \ref{def2}. We also show our main result in Theorem \ref{th-1} by using the assumptions on the interaction potential and relative entropy estimates. Here, we follow the blueprint of \cite{CL-2} being the most novel aspects how to deal with the case $m=1$ and the interaction potential. Finally, the last section is to remind the reader of the existence of weak solutions satisfying the needed properties for Theorem \ref{th-1} under suitable assumptions on the confinement potential. This part relies heavily on previous results in \cite{JAC17} being the most novel aspect how to deal with the confinement potential term.


\section{Relaxation limit: Relative entropy \& Convergence}

In this part, we devote ourselves to compare a  weak solution $(\rho,\rho\bu)$ of \eqref{Eqs} with a smooth solution $(\bar\rho,\bar{\rho}\bar{\bu})$ of \eqref{Eqs-bar-1} by using a relative entropy method.
The main result of this subsection is the following:
\begin{prop}\label{prop1}
Let $\Om$ be any smooth, connected, open subset of $\R^d$.
Let $(\rho,\rho\bu)$ be a  weak solution of \eqref{Eqs}  as in Definition \ref{def2} and $(\bar{\rho},\bar{\rho}\bar{\bu})$ be a smooth solution of \eqref{Eqs-bar-1}. Then
\begin{align}\label{diff-eqs-7}
&\int_{\Om}\left(\f{1}{\epz}h(\rho|\bar{\rho})  +  \f{1}{2}\rho|\bu-\bar{\bu}|^2  +  \f{C_k}{2\epz}(K\ast(\rho-\bar{\rho}))(\rho-\bar{\rho})\right)
\Big|_{\tau=0}^{\tau=t}dx   \nn  \\
&=-\f{1}{\epz}\int_0^t\int_{\Om}\rho|\bu-\bar{\bu}|^2dxd\tau     -   \int_0^t\int_{\Om}\rho\na_x\bar{\bu}:(\bu-\bar{\bu})\otimes(\bu-\bar{\bu})dxd\tau \nn \\
& \quad   -  \int_0^t\int_{\Om}\f{\rho}{\bar{\rho}}\bar{e}\cd(\bu-\bar{\bu})dxd\tau
     -  \f{1}{\epz}\int_0^t\int_{\Om}p(\rho|\bar{\rho})\na_x\cd\bar{\bu}dxd\tau   \nn \\
&  \quad - \f{C_k}{\epz}\int_0^t\int_{\Om}\left(K\ast(\rho-\bar{\rho})\right)\na_x\cd((\rho-\bar{\rho})\bar{\bu})  dxd\tau.
\end{align}
\end{prop}

 \textbf{Proof.} Firstly, we introduce the standard choice of test function in \eqref{weak-en-1}
\begin{equation}\label{theta}
\theta(\tau):=
\left\{
\begin{split}
&1,  \qquad \textrm{for} \,\, 0\le\tau<t, \\
&\f{t-\tau}{\kappa} + 1, \qquad \textrm{for}\,\, t\le \tau<t+\kappa,\\
&0,\qquad  \textrm{for}\,\,\tau\geq t+\kappa,
\end{split}
\right.
\end{equation}
and we have
\begin{align*}
  &\int_t^{t+\kappa}\int_{\Om}\f{1}{\kappa}\left(\f{1}{\epz}h(\rho) + \f12\rho|\bu|^2
   + \f{C_k}{2\epz}(K\ast\rho)\rho + \f{1}{\epz}\rho\Phi\right)dxd\tau   \nn \\
& \quad +   \f{1}{\epz}\int_0^t\int_{\Om}\rho|\bu|^2dxd\tau
     +   \f{1}{\epz}\int_t^{t+\kappa}\int_{\Om}\left(\f{t-\tau}{\kappa}+1\right)\rho|\bu|^2dxd\tau   \nn \\
  & \quad \quad=  \int_{\Om}\left(\f{1}{\epz}h(\rho) + \f12\rho|\bu|^2 + \f{C_k}{2\epz}(K\ast\rho)\rho + \f{1}{\epz}\rho\Phi\right)\Big|_{\tau=0}dx.
\end{align*}
Letting $\kappa$ tends to $0^{+}$, one has
\begin{align}\label{theta2}
\int_{\Om}\left(\f{1}{\epz}h(\rho) + \f12\rho|\bu|^2
   + \f{C_k}{2\epz}(K\ast\rho)\rho + \f{1}{\epz}\rho\Phi\right)\Big|_{\tau=0}^{\tau=t}dx
   =  -   \f{1}{\epz}\int_0^t\int_{\Om}\rho|\bu|^2dxd\tau.
\end{align}
Moreover, integrating \eqref{energy-2} over time interval $[0,t]$, one obtains
\begin{align}\label{theta3}
\int_{\Om}\left(\f{1}{\epz}h(\bar{\rho}) + \f12\bar{\rho}|\bar{\bu}|^2
   + \f{C_k}{2\epz}(K\ast\bar{\rho})\bar{\rho} + \f{1}{\epz}\bar{\rho}\Phi\right)\Big|_{\tau=0}^{\tau=t}dx
   =  -   \f{1}{\epz}\int_0^t\int_{\Om}\bar{\rho}|\bar{\bu}|^2dxd\tau
      +   \int_0^t\int_{\Om}\bar{\bu}\cd\bar{e}\,dxd\tau.
\end{align}

Next, we deduce from systems \eqref{Eqs} and \eqref{Eqs-bar-1} that the differences $\rho-\bar{\rho}$ and $\rho\bu-\bar{\rho}\bar{\bu}$ are given by the following equations
\begin{align}\label{Eqs-diff-1}
&\quad\p_t({\rho-\bar{\rho}})   +   \di_x(\rho\bu-\bar{\rho}\bar{\bu})  =  0,  \nn \\
&\p_t(\rho\bu-\bar{\rho}\bar{\bu})  +   \di_x(\rho\bu\otimes\bu - \bar{\rho}\bar{\bu}\otimes\bar{\bu})
    +  \f{1}{\epz}\na_x\left(p(\rho)-p(\bar{\rho})\right)   \nn \\
&\quad = -\f{C_k}{\epz}\left((\na_x K\ast\rho)\rho - (\na_x K\ast\bar{\rho})\bar{\rho}\right)   -   \f{1}{\epz}(\rho\bu-\bar{\rho}\bar{\bu})
   -   \f{1}{\epz}(\rho-\bar{\rho})\na_x\Phi     - \bar{e}.
\end{align}
Thus, the weak formulation for the equations satisfied by the differences $\rho-\bar{\rho}$ and $\rho\bu-\bar{\rho}\bar{\bu}$ in \eqref{Eqs-diff-1} reads
\begin{align}\label{Eqs-weak-1}
-\int_0^{\infty}\int_{\Om}\varphi_t(\rho-\bar{\rho})dxdt  -  \int_0^{\infty}\int_{\Om}\na_x\varphi\cd(\rho\bu-\bar{\rho}\bar{\bu})dxdt
    -  \int_{\Om}\varphi(\rho-\bar{\rho})\Big|_{t=0}dx    =    0,
\end{align}
\begin{align}\label{Eqs-weak-2}
&-\int_0^{\infty}\int_{\Om}\tilde{\varphi}_t\cd(\rho\bu-\bar{\rho}\bar{\bu})dxdt
-  \int_0^{\infty}\int_{\Om}\na_x\tilde{\varphi}:(\rho\bu\otimes\bu - \bar{\rho}\bar{\bu}\otimes\bar{\bu})dxdt    \nn \\
& - \f{1}{\epz} \int_0^{\infty}\int_{\Om} \di_x\tilde{\varphi}(p(\rho)-p(\bar{\rho})) dxdt
-  \int_{\Om}\tilde{\varphi}\cd(\rho\bu-\bar{\rho}\bar{\bu})\Big|_{t=0}dx      \nn \\
&  = -\f{C_k}{\epz}\int_0^{\infty}\int_{\Om}\tilde\varphi\cd\left((\na_x K\ast\rho)\rho - (\na_x K\ast\bar{\rho})\bar{\rho}\right)dxdt
-  \f{1}{\epz}\int_0^{\infty}\int_{\Om}\tilde{\varphi}\cd(\rho\bu - \bar{\rho}\bar{\bu})dxdt    \nn \\
&  \quad \, -  \f{1}{\epz}\int_0^{\infty}\int_{\Om}\tilde{\varphi}\cd  (\rho-\bar{\rho})\na_x\Phi dxdt
-  \int_0^{\infty}\int_{\Om}\tilde{\varphi}\cd  \bar{e}\,dxdt,
\end{align}
where $\varphi$ and $\tilde{\varphi}$ are Lipschitz test functions compactly supported in $[0,\infty)$ in time and $\tilde{\varphi}\cd\nu=0$ on $\p\Om$ when $\Om\neq\R^d$. Using the definition of $\theta(\tau)$ in \eqref{theta}, we introduce the test functions in the above relations
\begin{align*}
\varphi=\theta(\tau)\left(\f{1}{\epz}h'(\bar{\rho}) - \f12|\bar{\bu}|^2 + \f{C_k}{\epz}(K\ast\bar{\rho}) + \f{1}{\epz}\Phi\right),
\quad \tilde{\varphi}=\theta(\tau)\bar{\bu}
\end{align*}
and then we have by letting $\kappa\rightarrow0^{+}$ after substituting $\varphi$, $\tilde{\varphi}$ into \eqref{Eqs-weak-1} and \eqref{Eqs-weak-2}
\begin{align}\label{theta4}
\int_{\Om}&\left(\f{1}{\epz}h'(\bar{\rho}) - \f12|\bar{\bu}|^2 + \f{C_k}{\epz}(K\ast\bar{\rho}) + \f{1}{\epz}\Phi\right)
    (\rho - \bar{\rho})\Big|_{\tau=0}^{\tau=t}dx    \nn \\
&-\int_0^t\int_{\Om}\p_{\tau}\left(\f{1}{\epz}h'(\bar{\rho}) - \f12|\bar{\bu}|^2 + \f{C_k}{\epz}(K\ast\bar{\rho}) + \f{1}{\epz}\Phi\right)
    (\rho - \bar{\rho})dxd\tau   \nn \\
& - \int_0^t\int_{\Om}\na_x\left(\f{1}{\epz}h'(\bar{\rho}) - \f12|\bar{\bu}|^2 + \f{C_k}{\epz}(K\ast\bar{\rho}) + \f{1}{\epz}\Phi\right)
    \cd(\rho\bu - \bar{\rho}\bar{\bu})dxd\tau
=0
\end{align}
and
\begin{align}\label{theta5}
&\int_{\Om}\bar{\bu}\cd(\rho\bu - \bar{\rho}\bar{\bu})\Big|_{\tau=0}^{\tau=t}dx
- \int_0^t\int_{\Om}\p_{\tau}\bar{\bu}\cd(\rho\bu - \bar{\rho}\bar{\bu})dxd\tau     \nn \\
& -\int_0^t\int_{\Om}\na_x\bar{\bu}:(\rho\bu\otimes\bu - \bar{\rho}\bar{\bu}\otimes\bar{\bu})dxd\tau
-\f{1}{\epz}\int_0^t\int_{\Om}\di_x\bar{\bu}(p(\rho) - p(\bar{\rho}))dxd\tau   \nn \\
&=  - \f{C_k}{\epz}\int_0^t\int_{\Om}\bar{\bu}\cd\left((\na_x K\ast\rho)\rho - (\na_x K\ast\bar{\rho})\bar{\rho}\right)dxd\tau
   -  \f{1}{\epz}\int_0^t\int_{\Om}\bar{\bu}\cd(\rho\bu - \bar{\rho}\bar{\bu})dxd\tau   \nn \\
&\quad\, -  \f{1}{\epz}\int_0^t\int_{\Om}(\rho - \bar{\rho})\bar{\bu}\cd\na_x\Phi dxd\tau
  -  \int_0^t\int_{\Om}\bar{\bu}\cd\bar{e}\,dxd\tau.
\end{align}
We can deduce from the computation \eqref{theta2} $-$ \eqref{theta3} $-$ (\eqref{theta4} $+$ \eqref{theta5}) that
\begin{align}\label{theta6}
&\int_{\Om}\left(\f{1}{\epz}h(\rho|\bar{\rho})  +  \f{1}{2}\rho|\bu-\bar{\bu}|^2  +  \f{C_k}{2\epz}(K\ast(\rho-\bar{\rho}))(\rho-\bar{\rho})\right)
\Big|_{\tau=0}^{\tau=t}dx   \nn  \\
& =  -   \f{1}{\epz}\int_0^t\int_{\Om}\left(\rho|\bu|^2 - \bar{\rho}|\bar{\bu}|^2  -  \bar{\bu}\cd(\rho\bu-\bar{\rho}\bar{\bu})\right)dxd\tau   \nn \\
&\quad\, -\int_0^t\int_{\Om}\p_{\tau}\left(\f{1}{\epz}h'(\bar{\rho}) - \f12|\bar{\bu}|^2 + \f{C_k}{\epz}(K\ast\bar{\rho}) + \f{1}{\epz}\Phi\right)
    (\rho - \bar{\rho})dxd\tau
    - \int_0^t\int_{\Om}\p_{\tau}\bar{\bu}\cd(\rho\bu - \bar{\rho}\bar{\bu})dxd\tau   \nn \\
&\quad\, - \int_0^t\int_{\Om}\na_x\left(\f{1}{\epz}h'(\bar{\rho}) - \f12|\bar{\bu}|^2 + \f{C_k}{\epz}(K\ast\bar{\rho}) + \f{1}{\epz}\Phi\right)
    \cd(\rho\bu - \bar{\rho}\bar{\bu})dxd\tau  \nn \\
&\quad\, -\int_0^t\int_{\Om}\na_x\bar{\bu}:(\rho\bu\otimes\bu - \bar{\rho}\bar{\bu}\otimes\bar{\bu})dxd\tau
-\f{1}{\epz}\int_0^t\int_{\Om}\di_x\bar{\bu}(p(\rho) - p(\bar{\rho}))dxd\tau   \nn \\
&\quad\, + \f{C_k}{\epz}\int_0^t\int_{\Om}\bar{\bu}\cd\left((\na_x K\ast\rho)\rho - (\na_x K\ast\bar{\rho})\bar{\rho}\right)dxd\tau
 +  \f{1}{\epz}\int_0^t\int_{\Om}(\rho - \bar{\rho})\bar{\bu}\cd\na_x\Phi dxd\tau.
\end{align}
Deducing from \eqref{Eqs-bar-1} by using $\bar{\rho}>0$, one can obtain the equation satisfied by $\bar{\bu}$
\begin{align}\label{Eqs-bar-u}
\p_{\tau}\bar{\bu}   +   \bar{\bu}\cd\na_x\bar{\bu}
=- \f{1}{\epz}\na_x h'(\bar{\rho})  -  \f{C_k}{\epz}\na_x(K\ast\bar{\rho})  -  \f{1}{\epz}\bar{\bu}  -  \f{1}{\epz}\na_x\Phi  +  \f{\bar{e}}{\bar{\rho}},
\end{align}
where we have used \eqref{h-p-1}. Furthermore, multiplying \eqref{Eqs-bar-u} with $\rho(\bu-\bar{\bu})$ leads to
\begin{align}\label{Eqs-bar-u1}
&\p_{\tau}\big(-\f12|\bar{\bu}|^2\big)(\rho-\bar{\rho})   +   \p_{\tau}\bar{\bu}\cd(\rho\bu-\bar{\rho}\bar{\bu})
+  \na_x\big(-\f12|\bar{\bu}|^2\big)\cd(\rho\bu-\bar{\rho}\bar{\bu})
 + \na_x\bar{\bu}:(\rho\bu\otimes\bu - \bar{\rho}\bar{\bu}\otimes\bar{\bu})   \nn \\
&\quad = \rho\na_x\bar{\bu}:(\bu-\bar{\bu})\otimes(\bu-\bar{\bu})  -  \f{1}{\epz}\rho\na_x h'(\bar{\rho})\cd(\bu-\bar{\bu})
  -  \f{C_k}{\epz}\rho\na_x(K\ast\bar{\rho})\cd(\bu-\bar{\bu})    \nn \\
&\qquad    -  \f{1}{\epz}\rho\bar{\bu}\cd(\bu-\bar{\bu})
-  \f{1}{\epz}\rho\na_x\Phi\cd(\bu-\bar{\bu})  +\f{\rho}{\bar{\rho}}\bar{e}\cd(\bu-\bar{\bu}).
\end{align}
Substituting \eqref{Eqs-bar-u1} into \eqref{theta6} and using \eqref{Eqs-bar-1}$_1$, one gets
\begin{align}\label{theta7}
&\int_{\Om}\left(\f{1}{\epz}h(\rho|\bar{\rho})  +  \f{1}{2}\rho|\bu-\bar{\bu}|^2  +  \f{C_k}{2\epz}(K\ast(\rho-\bar{\rho}))(\rho-\bar{\rho})\right)
\Big|_{\tau=0}^{\tau=t}dx   \nn  \\
&=-\f{1}{\epz}\int_0^t\int_{\Om}\rho|\bu-\bar{\bu}|^2dxd\tau     -   \int_0^t\int_{\Om}\rho\na_x\bar{\bu}:(\bu-\bar{\bu})\otimes(\bu-\bar{\bu})dxd\tau \nn \\
&\quad + \f{C_k}{\epz}\int_0^t\int_{\Om}\left(\na_x K\ast(\rho-\bar{\rho})\right)\cd\rho\bar{\bu}dxd\tau
   -  \int_0^t\int_{\Om}\f{\rho}{\bar{\rho}}\bar{e}\cd(\bu-\bar{\bu})dxd\tau     \nn \\
&\quad -  \f{1}{\epz}\int_0^t\int_{\Om}p(\rho|\bar{\rho})\di_x\bar{\bu}dxd\tau
   -   \int_0^t\int_{\Om}(\rho-\bar{\rho})\p_{\tau}\left(\f{C_k}{\epz}(K\ast\bar{\rho})  +  \f{1}{\epz}\Phi\right)dxd\tau.
\end{align}
Due to the fact that $K$ is symmetric, one can deduce that
\begin{align*}
\int_{\Om}(K\ast\rho)\bar{\rho}dx  =  \int_{\Om}(K\ast\bar{\rho})\rho dx,
\end{align*}
consequently,
\begin{align}\label{theta8}
&-   \int_0^t\int_{\Om}(\rho-\bar{\rho})\p_{\tau}\left(\f{C_k}{\epz}(K\ast\bar{\rho})  +  \f{1}{\epz}\Phi\right)dxd\tau   \nn \\
&=  -   \f{C_k}{\epz}\int_0^t\int_{\Om}(\rho-\bar{\rho})\p_{\tau}(K\ast\bar{\rho})  dxd\tau
=  -   \f{C_k}{\epz}\int_0^t\int_{\Om}\left(K\ast(\rho-\bar{\rho})\right)\p_{\tau}\bar{\rho}  dxd\tau    \nn \\
&=   \f{C_k}{\epz}\int_0^t\int_{\Om}\left(K\ast(\rho-\bar{\rho})\right)\di_x(\bar{\rho}\bar{\bu})  dxd\tau    \nn \\
&=   \f{C_k}{\epz}\int_0^t\int_{\Om}\left(K\ast(\rho-\bar{\rho})\right)\di_x(\rho\bar{\bu})  dxd\tau
   - \f{C_k}{\epz}\int_0^t\int_{\Om}\left(K\ast(\rho-\bar{\rho})\right)\di_x((\rho-\bar{\rho})\bar{\bu})  dxd\tau.
\end{align}
Hence, one can finally obtain by substituting \eqref{theta8} into \eqref{theta7} that
\begin{align*}
&\int_{\Om}\left(\f{1}{\epz}h(\rho|\bar{\rho})  +  \f{1}{2}\rho|\bu-\bar{\bu}|^2  +  \f{C_k}{2\epz}(K\ast(\rho-\bar{\rho}))(\rho-\bar{\rho})\right)
\Big|_{\tau=0}^{\tau=t}dx   \nn  \\
&=-\f{1}{\epz}\int_0^t\int_{\Om}\rho|\bu-\bar{\bu}|^2dxd\tau     -   \int_0^t\int_{\Om}\rho\na_x\bar{\bu}:(\bu-\bar{\bu})\otimes(\bu-\bar{\bu})dxd\tau \nn \\
&\quad  \,  -  \int_0^t\int_{\Om}\f{\rho}{\bar{\rho}}\bar{e}\cd(\bu-\bar{\bu})dxd\tau
     -  \f{1}{\epz}\int_0^t\int_{\Om}p(\rho|\bar{\rho})\di_x\bar{\bu}dxd\tau   \nn \\
&\quad \,  - \f{C_k}{\epz}\int_0^t\int_{\Om}\left(K\ast(\rho-\bar{\rho})\right)\di_x((\rho-\bar{\rho})\bar{\bu})  dxd\tau.
\end{align*}
This exactly completes the proof of the Proposition \ref{prop1}.        \qquad $\Box$


\subsection{Convergence in the relaxation limit}
In this subsection, we will establish the convergence property in the relaxation limit from \eqref{Eqs} to \eqref{Eqs-bar-1} based on Proposition \ref{prop1}.

With the relative relation \eqref{diff-eqs-7} between solutions to \eqref{Eqs} and \eqref{Eqs-bar-1} at hand, we can prove Theorem \ref{th-1} by showing that terms on the right-hand-side of \eqref{diff-eqs-7} can be absorbed or are $O(\ep)$.

Before getting into the proof of our main theorem, we need firstly to have some auxiliary lemmas which essentially indicate that the relative potential energy can be bounded from below by some positive functions.
\begin{lem}\label{lem-2}
Let $h(\rho)$ be defined by \eqref{h-p-1} and \eqref{h-1}. Then for any $\bar{\rho}>0$, we have  the following estimates
\begin{align}\label{h-2}
h(\rho|\bar{\rho})
\ge\f{k_1}{2}\min\left\{\f{1}{\rho},\f{1}{\bar{\rho}}\right\}|\rho-\bar{\rho}|^2 \quad \,\,for\,\, any\,\,0<\rho<\infty\,\,and\,\, m=1
\end{align}
and
\begin{align}\label{h-3}
h(\rho|\bar{\rho})
\ge \f{k_2m}{2}\min\{\rho^{m-2},\bar{\rho}^{m-2}\}|\rho-\bar{\rho}|^2 \quad \,\,for\,\, any\,\,0<\rho<\infty\,\,and\,\, 1<m\le2.
\end{align}
\end{lem}
\textbf{Proof.} For the case of $m=1$, the Taylor expansion of $h(\rho)$ at $\bar{\rho}$ reads
\begin{align*}
h(\rho)
=h(\bar{\rho})  +   h'(\bar{\rho})(\rho-\bar{\rho})   +   \f{h''(\rho_{*})}{2}|\rho-\bar{\rho}|^2,
   \quad \rho_*\in [\rho,\bar{\rho}],
\end{align*}
which implies
\begin{align*}
h(\rho|\bar{\rho})
=\f{h''(\rho_{*})}{2}|\rho-\bar{\rho}|^2
=\f{k_1}{2\rho_{*}}|\rho-\bar{\rho}|^2
\ge\f{k_1}{2}\min\left\{\f{1}{\rho},\f{1}{\bar{\rho}}\right\}|\rho-\bar{\rho}|^2.
\end{align*}
For the case of $1<m\le2$, similarly, the Taylor expansion of $h(\rho)$ at $\bar{\rho}$ entails that
\begin{align*}
h(\rho|\bar{\rho})
=\f{h''(\xi)}{2}|\rho-\bar{\rho}|^2
=\f{k_2m}{2}\xi^{m-2}|\rho-\bar{\rho}|^2
\ge \f{k_2m}{2}\min\{\rho^{m-2},\bar{\rho}^{m-2}\}|\rho-\bar{\rho}|^2  \quad (\xi\in [\rho,\bar{\rho}]).
\end{align*}
This completes the proof of \eqref{h-2} and \eqref{h-3}.        \qquad $\Box$
\\

We remind the readers a result proved in \cite[Lemma 2.4]{CL-1}.
\begin{lem}\label{lem-1}
Let $h(\rho)$ be defined by \eqref{h-p-1} and \eqref{h-1}. If $\bar{\rho}\in I=[\underline{\del},\overline{\del}]$ with $\underline{\del}>0$ and $\overline{\del}<+\infty$, $m>1$, then there exist positive constants $R_0$ (depending on $I$) and $C_1$, $C_2$ (depending on $I$ and $R_0$) such that
\begin{equation*}
h(\rho|\bar{\rho})\geq
\left\{
\begin{split}
&C_1 |\rho-\bar{\rho}|^2  \qquad for \,\, 0\le\rho\le R_0, \,\bar{\rho}\in I,  \\
&C_2 |\rho-\bar{\rho}|^m  \qquad for \,\, \rho> R_0, \,\bar{\rho}\in I, m>1.
\end{split}
\right.
\end{equation*}
\end{lem}

Given $h(\rho)$ defined by \eqref{h-p-1} and \eqref{h-1},
we can verify by using a similar way as in \cite[Lemma 2.3]{CL-1} that
\begin{align}\label{re-p}
|p(\rho|\bar{\rho})|   \le   C h(\rho|\bar{\rho})   \qquad \forall \rho,\bar{\rho}>0, \mbox{ and for some }C>0.
\end{align}

\begin{lem}\label{lem4}
Let $\Om$ be any smooth, connected, open subset of $\R^d$ and $h(\rho)$ be defined by \eqref{h-p-1} and \eqref{h-1}.
Assume one of the following conditions hold:\\
(i) If $2-\f{2}{d}\le m \le 2$ with $d\ge2$
and the interaction potential $K$ satisfies
$K\in L^{\f{m}{2(m-1)}}(\Om)\cap W^{1,\infty}(\Om)$,\\
\\
(ii) If $\Om=\mathbb{T}^d$ or $\Om$ is a bounded domain in $\R^d$, $m\ge2-\f{2}{d}$ with $d\ge2$, $\bar{\rho}\in[\underline{\del},\overline{\del}]$ with $\underline{\del}>0$ and $\overline{\del}<\infty$ and $K$ satisfies
$K\in L^p(\Om) \quad (1<p\le\infty)$.\\
\\
Then there exists a positive constant $C_{*}$ such that
\begin{align}\label{energy-5}
\left|\int_{\Om}(\rho-\bar{\rho})(K\ast(\rho-\bar{\rho}))dx\right|
\le C_{*}\int_{\Om}h(\rho|\bar{\rho})dx \quad \mbox{ for a.a.}t\in[0,T].
\end{align}
\end{lem}
\textbf{Proof.} Firstly, let us work with the case $m=1$ and $d=2$. By using H\"{o}lder's inequality and Young's inequality, we obtain
\begin{align}\label{en-1}
\left|\int_{\Om}(\rho-\bar{\rho})(K\ast(\rho-\bar{\rho}))dx\right|
\le C\|K\|_{L^{\infty}(\Om)}\|\rho-\bar{\rho}\|^2_{L^1(\Om)}.
\end{align}
Due to
\begin{align}\label{en-2}
\|\rho-\bar{\rho}\|_{L^1(\Om)}
=\int_{\Om}|\rho-\bar{\rho}|dx
&=\int_{\Om}\sqrt{\min\left\{\f{1}{\rho},\f{1}{\bar{\rho}}\right\}} \, |\rho-\bar{\rho}| \,
     \left(\sqrt{\min\left\{\f{1}{\rho},\f{1}{\bar{\rho}}\right\}}\right)^{-1}dx \nn \\
&\le \left(\int_{\Om}\min\left\{\f{1}{\rho},\f{1}{\bar{\rho}}\right\}|\rho-\bar{\rho}|^2dx\right)^{\f12}
    \left(\int_{\Om}\max\{\rho,\bar{\rho}\}dx\right)^{\f12}  \nn \\
&\le C\left(\int_{\Om}\min\left\{\f{1}{\rho},\f{1}{\bar{\rho}}\right\}|\rho-\bar{\rho}|^2dx\right)^{\f12},
\end{align}
where we have used the mass conservation property of $\rho$ and $\bar{\rho}$ in the last inequality. We can claim by substituting \eqref{en-2}
into \eqref{en-1} and using \eqref{h-2} that \eqref{energy-5} is valid for $m=1$, $d=2$.

For the case of $1< m\le2$ with $d=2$ and $2-\f{2}{d}\le m\le2$ with $d\ge3$,   we have
\begin{align}\label{en-3}
\left|\int_{\Om}(\rho-\bar{\rho})(K\ast(\rho-\bar{\rho}))dx\right|
\le C\|K\|_{L^{\f{m}{2(m-1)}}(\Om)}\|\rho-\bar{\rho}\|^2_{L^m(\Om)}.
\end{align}
Since $\Phi$ is bounded from blow and $\int_{\Om}(K\ast\rho)\rho\, dx \le \|K\|_{L^{\infty}(\Om)}\|\rho\|^2_{L^1(\Om)}$, one can deduce from the energy estimates \eqref{energy-1} and \eqref{energy-2} that $\int_{\Om}\rho^mdx$ and $\int_{\Om}\bar{\rho}^mdx$ are bounded.
Thus we have
\begin{align*}
\|\rho-\bar{\rho}\|^m_{L^m(\Om)}
&=\int_{\Om}|\rho-\bar{\rho}|^mdx   \nn \\
&=\int_{\Om}\left(\f{k_2m}{2}\min\{\rho^{m-2},\bar{\rho}^{m-2}\}\right)^{\f{m}{2}}|\rho-\bar{\rho}|^m
   \left(\f{k_2m}{2}\min\{\rho^{m-2},\bar{\rho}^{m-2}\}\right)^{-\f{m}{2}}dx   \nn \\
&\le  \left(\f{k_2m}{2}\right)^{-\f{m}{2}}\left(\int_{\Om}\f{k_2m}{2}\min\{\rho^{m-2},\bar{\rho}^{m-2}\}  |\rho-\bar{\rho}|^2dx\right)^{\f{m}{2}}
     \left(\int_{\Om}\max\{\rho^m,\bar{\rho}^m\}dx\right)^{\f{2-m}{2}}   \nn \\
&\le C\left(\int_{\Om}\f{k_2m}{2}\min\{\rho^{m-2},\bar{\rho}^{m-2}\}  |\rho-\bar{\rho}|^2dx\right)^{\f{m}{2}},
\end{align*}
which implies that
\begin{align}\label{J-3-6}
\|\rho-\bar{\rho}\|^2_{L^m(\Om)}
\le C\int_{\Om}\f{k_2m}{2}\min\{\rho^{m-2},\bar{\rho}^{m-2}\}  |\rho-\bar{\rho}|^2dx.
\end{align}
Substituting \eqref{J-3-6} into \eqref{en-3} and using \eqref{h-3}, then, for $1<m\le2$ with $d=2$ and $2-\f{2}{d}\le m\le2$ with $d\ge3$, the proof of \eqref{energy-5} is completed.

It remains to prove the case of $m>2$ with any $d\ge 2$ when $\Om=\mathbb{T}^d$ or $\Om$ is a bounded domain. In Lemma \ref{lem-1}, by enlarging if necessary $R_0$ so that $|\rho-\bar{\rho}|\ge1$ for $\rho>R_0$ and $\bar{\rho}\in[\underline{\delta},\overline{\del}]$, then we have
\begin{align*}
h(\rho|\bar{\rho})\ge C|\rho-\bar{\rho}|^2,\quad \textrm{for}\,\,m>2,\,\,\rho\ge0,\,\,\bar{\rho}\in[\underline{\delta},\overline{\del}].
\end{align*}
Thus, one deduce that
\begin{align*}
\left|\int_{\Om}(\rho-\bar{\rho})(K\ast(\rho-\bar{\rho}))dx\right|
\le C\|K\|_{L^{\f{r}{2(r-1)}}(\Om)}\|\rho-\bar{\rho}\|^2_{L^r(\Om)}
\le C\|\rho-\bar{\rho}\|^2_{L^2(\Om)}
\le C\int_{\Om}h(\rho|\bar{\rho})dx,
\end{align*}
where $1\le r<2$ and we have used the fact that $\Om=\mathbb{T}^d$ or $\Om$ is a bounded domain in the last second inequality.
The proof of \eqref{energy-5} is completed.       \qquad $\Box$

\begin{coro}\label{coro1}
Let the assumptions in Lemma \ref{lem4} hold and the parameter $C_k$ is such that
$C_{k}<\f{2}{C_{*}}$,
where $C_{*}$ is defined in \eqref{energy-5}, then for $\lambda:=1-\f{C_kC_{*}}{2}>0$
\begin{align*}
\int_{\Om}h(\rho|\bar{\rho})   +   \f{C_k}{2}\int_{\Om}(\rho-\bar{\rho})(K\ast(\rho-\bar{\rho}))dx
\ge \lambda \int_{\Om}h(\rho|\bar{\rho}) \quad \mbox{ for a.a.}t\in[0,T].
\end{align*}
\end{coro}

So far, all the preparations have been done, we now start to prove our main result.\\
\\
\textbf{Proof of Theorem \ref{th-1}.} Firstly, one can easily see from the definition of $\Theta(t)$ in \eqref{varphi-1} and the relative entropy identity \eqref{diff-eqs-7} that
\begin{align}\label{varphi-2}
\Theta(t)   +   \f{1}{\epz}\int_0^t\int_{\Om}\rho|\bu-\bar{\bu}|^2dxd\tau
&=  \Theta(0)   -   \int_0^t\int_{\Om}\rho\na_x\bar{\bu}:(\bu-\bar{\bu})\otimes(\bu-\bar{\bu})dxd\tau \nn \\
&\quad   - \f{C_k}{\epz}\int_0^t\int_{\Om}\left(K\ast(\rho-\bar{\rho})\right)\di_x((\rho-\bar{\rho})\bar{\bu})  dxd\tau   \nn \\
&\quad     -  \f{1}{\epz}\int_0^t\int_{\Om}p(\rho|\bar{\rho})\di_x\bar{\bu}dxd\tau
 -  \int_0^t\int_{\Om}\f{\rho}{\bar{\rho}}\bar{e}\cd(\bu-\bar{\bu})dxd\tau    \nn \\
&:= \Theta(0) +  J_1  +   J_2   +  J_3  +  J_4.
\end{align}
Now, we estimate $J_1$, $J_2$, $J_3$, and $J_4$ one by one. Using the relation between $p$ and $h$ in \eqref{h-p-1} and the definition of $\bar{\bu}$ in \eqref{bar-m}, then we deduce that
$\bar{\bu}=-\na_x h'(\bar{\rho})   - C_k(\na_x K\ast\bar{\rho})  -  \na_x\Phi$
and
$\na_x\bar{\bu}$ are bounded functions due to the smoothness assumption on $\bar{\rho}$.

For $J_1$, one obtains
\begin{align}\label{J-1}
J_1
&=-\int_0^t\int_{\Om}\rho\na_x\bar{\bu}:(\bu-\bar{\bu})\otimes(\bu-\bar{\bu}) dxd\tau   \nn \\
&\le   \|\na_x\bar{\bu}\|_{L^{\infty}((0,T)\times\Om)}\int_0^t\int_{\Om}\rho|\bu-\bar{\bu}|^2dxd\tau
\le  C\int_0^t\Theta(\tau)d\tau.
\end{align}

We will estimate $J_2$ for three different cases. The first case is for $m=1$ and $d=2$, the second case is for $2-\f{2}{d}<m\le2$ with $d\ge2$ or $2-\f{2}{d}\le m\le 2$ with $d\ge 3$ and the third case is for $m>2$ for any $d\ge 2$.  For $m=1$ and $d=2$,   using H\"{o}lder's inequality and Young's inequality, one deduces by using integration by parts that
\begin{align}\label{J-3-1}
J_2
&=- \f{C_k}{\epz}\int_0^t\int_{\Om}\di_x((\rho-\bar{\rho})\bar{\bu}) \left(K\ast(\rho-\bar{\rho})\right)dxd\tau
=  \f{C_k}{\epz}\int_0^t\int_{\Om}(\rho-\bar{\rho})\bar{\bu}\cd \left(\na_x K\ast(\rho-\bar{\rho})\right)dxd\tau   \nn \\
&\le \f{C}{\epz}\int_0^t\|\bar{\bu}\|_{L^{\infty}(\Om)}  \|\na_x K\|_{L^{\infty}(\Om)}  \|\rho-\bar{\rho}\|^2_{L^1(\Om)}d\tau
\le \f{C}{\epz}\int_0^t\|\rho-\bar{\rho}\|^2_{L^1(\Om)}d\tau    \nn \\
&\le \f{C}{\epz}\int_0^t\int_{\Om}\min\left\{\f{1}{\rho},\f{1}{\bar{\rho}}\right\}|\rho-\bar{\rho}|^2dxd\tau
\le \f{C}{\epz}\int_0^t\int_{\Om} h(\rho|\bar{\rho})d\tau
\le C\int_0^t\Theta(\tau)d\tau,
\end{align}
where we have used \eqref{en-2} in the last third inequality and Lemma \ref{lem-2} in the last second inequality.

For the case $1<m\le2$ with $d=2$ and $2-\f{2}{d}\le m\le 2$ with $d\ge 3$,  }  we obtain by using interpolation inequality that
\begin{align}\label{J-3-4}
J_2
&=  \f{C_k}{\epz}\int_0^t\int_{\Om}(\rho-\bar{\rho})\bar{\bu}\cd \left(\na_x K\ast(\rho-\bar{\rho})\right)dxd\tau  \nn \\
&\le \f{C_k}{\epz}\|\bar{\bu}\|_{L^{\infty}(0,T;L^{\f{m}{2(m-1)}}(\Om))}    \|\na_x K\|_{L^{\infty}(\Om)}
       \int_0^t\|\rho-\bar{\rho}\|^2_{L^m(\Om)}d\tau
\le \f{C}{\epz}\int_0^t \|\rho-\bar{\rho}\|^2_{L^m(\Om)}d\tau.
\end{align}
Substituting \eqref{J-3-6} into \eqref{J-3-4}, we have by Lemma \ref{lem-2}
\begin{align}\label{J-3-7}
J_2
\le\f{C}{\epz}\int_0^t\int_{\Om}\f{m}{2}\min\{\rho^{m-2},\bar{\rho}^{m-2}\}  |\rho-\bar{\rho}|^2dxd\tau
\le\f{C}{\epz}\int_0^t\int_{\Om}h(\rho|\bar{\rho})dxd\tau
\le C\int_0^t\Theta(\tau)d\tau.
\end{align}

Finally, for the case $m>2$ and any $d\ge 2$,
we have
\begin{align}\label{J-3-5}
J_2
&=  \f{C_k}{\epz}\int_0^t\int_{\Om}(\rho-\bar{\rho})\bar{\bu}\cd \left(\na_x K\ast(\rho-\bar{\rho})\right)dxd\tau   \nn \\
&\le \f{C_k}{\epz}\|\bar{\bu}\|_{L^{\infty}(0,T;L^{p}(\Om))}    \|\na_x K\|_{L^q(\Om)}
       \int_0^t\|\rho-\bar{\rho}\|^2_{L^{2}(\Om)}d\tau   \nn \\
&\le  \f{C}{\epz}\int_0^t\|\rho-\bar{\rho}\|^2_{L^2(\Om)}d\tau
\le \f{C}{\epz}\int_0^t\int_{\Om} h(\rho|\bar{\rho})d\tau
\le C\int_0^t\Theta(\tau)d\tau,
\end{align}
where $\f{1}{p}+\f{1}{q}=1$, due to Lemma \ref{lem-1} used in the last second inequality.

For $J_3$, by \eqref{re-p}, one has
\begin{align}\label{J-2}
J_3
&=-\f{1}{\epz}\int_0^t\int_{\Om}p(\rho|\bar{\rho})\di_x\bar{\bu} dxd\tau
\le \f{1}{\epz}\|\na_x\bar{\bu}\|_{L^{\infty}((0,T)\times\Om)}  \int_0^t\int_{\Om}|p(\rho|\bar{\rho})|dxd\tau   \nn \\
&\le  \f{C}{\epz}\int_0^t \int_{\Om}h(\rho|\bar{\rho})dxd\tau
\le C\int_0^t\Theta(\tau)d\tau.
\end{align}

For $J_4$, we similarly have that
\begin{align}\label{J-5}
J_4
&=-\int_0^t\int_{\Om}\rho(\bu-\bar{\bu})\cd\f{\bar{e}}{\bar{\rho}}dxd\tau
\le \f{1}{2\epz}\int_0^t\int_{\Om}\rho|\bu-\bar{\bu}|^2dxd\tau
     +    \f{\epz}{2}\int_0^t\int_{\Om}\rho\left|\f{\bar{e}}{\bar{\rho}}\right|^2dxd\tau   \nn \\
&\le \f{1}{2\epz}\int_0^t\int_{\Om}\rho|\bu-\bar{\bu}|^2dxd\tau    +    C\epz t,
\end{align}
where we have used the fact that $\bar{e}$ is bounded and the mass conservation of $\rho$ in the last inequality.
Substituting \eqref{J-1}, \eqref{J-3-1}, \eqref{J-3-7}, \eqref{J-3-5}, \eqref{J-2} and \eqref{J-5} into \eqref{varphi-2}, one can see that
\begin{equation*}
\Theta(t)  +   \f{1}{2\epz}\int_0^t\int_{\Om}\rho|\bu-\bar{\bu}|^2dxd\tau
\le\Theta(0)  +   C\int_0^t\Theta(\tau)d\tau   +  C\epz t.
\end{equation*}
Hence, Gronwall's inequality leads to
\begin{equation*}
\Theta(t)
\le
\tilde{C}(\Theta(0)  +  \epz)
\end{equation*}
for any $t\in(0,T]$, where $\tilde{C}$ is a positive constant depending on $T$. This completes the proof of Theorem \ref{th-1}.     \qquad $\Box$

Recalling the definition of $\Theta(t)$ in \eqref{varphi-1} and the properties of $h(\rho|\bar{\rho})$ showed in Lemma \ref{lem-2} and Lemma \ref{lem-1}, we can easily conclude the following result.

\begin{coro}\label{coro2}
Let all conditions in Theorem \ref{th-1} hold, then we can conclude that the weak solution  of \eqref{Eqs-1} converges to the solution $(\bar{\rho},\bar{\rho}\bar{\bu})$ of \eqref{equi-eq-1} in the sense that
\begin{align*}
\|\rho-\bar{\rho}\|_{L^{\infty}(0,T;L^2(\Om))}\rightarrow 0 \qquad as \quad\epz\rightarrow 0
\end{align*}
and
\begin{align*}
\|\sqrt{\rho}(\bu-\bar{\bu})\|_{L^{\infty}(0,T;L^2(\Om))\cap L^2(0,T;L^2(\Om))} \rightarrow 0 \qquad as \quad\epz\rightarrow 0,
\end{align*}
where $\bar{\bu}=-\na_x h(\bar{\rho})  -  C_k(\na_x K \ast\bar{\rho})  -  \na_x\Phi$.
\end{coro}


\section{Weak solutions to the Hydrodynamic system}

Our goal in this section is to prove existence of weak solutions to the system \eqref{Eqs} by using the methods of convex integration and oscillatory lemma shown in the seminal work by C. De Lellis and L. Sz\'{e}kelyhidi \cite{DL10}. Similar methods are later applied to deal with the compressible Euler system by E. Chiodaroli \cite{CH14}, the Euler systems with non-local interactions by J. A. Carrillo et al. \cite{JAC17} and some more general "variable coefficients" problems in \cite{DD15,CH15,EF15,EF16}.

The proof of the existence theory for the weak solutions of Euler flow \eqref{Eqs} on any bounded domain $\Om\subset\R^d$, $d=2,3$ with smooth boundary can be done by adapting the method of convex integration in \cite{JAC17}. Solvability for other cases mentioned in this paper, i.e. $\Om\subset\R^d$ ($d\ge2$) unbounded or $\Om\subset\R^d$ ($d\ge4$) bounded with smooth boundary, are left open.

For simplicity, we take the coefficients $\epz=C_k=1$ in \eqref{Eqs} and restrict ourselves to the spatially periodic boundary conditions, i.e. $x\in\Om$, where \begin{align}\label{domain-b}
\Om=([-1,1]|_{\{-1,1\}})^{d}, \quad d=2,3,
\end{align}
is the "flat" torus. One should notice that this method is applicable for the general connected bounded domains $\Om\subset\R^d$ with smooth boundary endowed with the no-flux boundary conditions $\bu\cd\nu|_{\p\Om}=0$.
Thus, we consider the solvability of the following system
\begin{equation}\label{Eqs-solv}
\begin{split}
\p_t\rho   +    \di_x (\rho\bu)   &=0  , \\
\p_t(\rho\bu)  +  \di_x(\rho\bu\otimes\bu) +  \na_x p(\rho)
&=  - (\na_x K\ast\rho)\rho  -  \rho\bu
         -  \rho\na_x\Phi
\end{split}
\end{equation}
with initial data
\begin{align}\label{initial-solv}
\rho(0,\cd)=\rho_0, \qquad \bu(0,\cd)=\bu_0.
\end{align}

\begin{thm}\label{thm2}
Let $T>0$ be given and $d=2,3$. Suppose that
\begin{align*}
p\in C[0,\infty)\cap C^2(0,\infty), \quad p(0)=0, \quad K\in C^2(\Om), \quad \Phi\in C^2(\Om).
\end{align*}
Let the initial data $\rho_0$, $\bu_0$ satisfy $\rho_0\in C^2(\Om)$, $\rho_0\ge\underline{\rho}>0$ in $\Om$, $\bu_0\in C^3(\Om;\R^d)$. Then the system
\eqref{Eqs-solv}, \eqref{initial-solv}, \eqref{domain-b} admits infinitely many solutions in the space-time cylinder $(0,T)\times\Om$ belonging to the class
\begin{align*}
\rho\in C^2([0,T]\times\Om),\quad \rho>0, \quad \bu\in C_{\textrm{weak}}([0,T];L^2(\Om;\R^d))\cap L^{\infty}((0,T)\times\Om;\R^d).
\end{align*}
\end{thm}

For the reader's convenience and completeness of this paper, we give a sketch of the proof of Theorem \ref{thm2} following the blueprint of \cite{JAC17}.


\subsection{Solvability of the abstract Euler system}
Firstly, we introduce the notations
\begin{align*}
\bv\otimes\bw\in\R_{\textrm{sym}}^{d\times d}, \quad [\bv\otimes\bw]_{i,j}=v_iv_j,\quad \textrm{and} \quad
  \bv\odot\bw\in\R_{\textrm{sym},0}^{d\times d}, \quad \bv\odot\bw = \bv\otimes\bw - \f{1}{d}\bv\cd\bw\mathbb{I},
\end{align*}
where $\bv,\bw\in\R^d$ are two vectors, $\R_{\textrm{sym}}^{d\times d}$ denotes the space of $d\times d$ symmetric matrices over the Euclidean space $\R^d$, $d=2,3$, $\R_{\textrm{sym,0}}^{d\times d}$ means its subspace of those with zero trace.
Recalling the abstract result in \cite{DL10,EF15} which will be used later to prove our existence result, we consider the following abstract Euler form:

Find a vector field $\bv\in C_{\textrm{weak}}([0,T];L^2(\Om;\R^d))$ satisfying
\begin{align}\label{abeqs1}
\p_t\bv  +  \di_x\left(\f{(\bv + \bh[\bv])\odot(\bv + \bh[\bv])}{r[\bv]}  +   \mathbb{H}[\bv]\right)   =   0,
\quad \di_x\bv = 0
\end{align}
in $\mathcal{D}'((0,T)\times\Om;\R^d)$,
\begin{align}\label{abeqs2}
\f12\f{|\bv + \bh[\bv]|^2}{r[\bv]}(t,x)  =   e[\bv](t,x) \quad \textrm{for} \,\, \textrm{a.a.}\,\, (t,x)\in (0,T)\times\Om,
\end{align}
\begin{align}\label{abeqs3}
\bv(0,\cd)=\bv_0, \quad \bv(T,\cd)=\bv_T,
\end{align}
where $\bh[\bv]$, $r[\bv]$, $\mathbb{H}[\bv]$, and $e[\bv]$ are given (nonlinear) operators.

\begin{dnt}\label{def3}
Let $Q\subset(0,T)\times\Om$ be an open set such that
\begin{align*}
|Q|=|(0,T)\times\Om|.
\end{align*}
An operator
\begin{align*}
b: C_{\textrm{weak}}([0,T];L^2(\Om;\R^d))\cap L^{\infty}((0,T)\times\Om;\R^d)  \rightarrow C_b(Q,\R^m)
\end{align*}
is $Q$-continuous if
\begin{itemize}
  \item b maps bounded sets in $L^{\infty}((0,T)\times\Om;\R^d)$ on bounded sets in $C_b(Q,\R^m)$;
  \item b is continuous, specifically,
  \begin{align*}
  b[\bv_n] \rightarrow b[\bv] \,\, in \,\, C_b(Q,\R^m) \,\,(uniformly\,\, for \,\, (t,x)\in Q)
  \end{align*}
  \begin{center}
  whenever
  \end{center}
  \begin{align*}
  \bv_n \rightarrow \bv \,\, in\,\, C_{weak}([0,T];L^2(\Om;\R^d))\,\, and \,\, weakly-(*)\,\, in\,\ L^{\infty}((0,T)\times\Om;\R^d);
  \end{align*}
  \item b is causal (non-anticipative), meaning
  \begin{align*}
  \bv(t,\cd)=\bw(t,\cd)\,\, for \,\, 0\le t\le \tau\le T\,\, implies\,\, b[\bv]=b[\bw]\,\, in \,\, [(0,\tau]\times\Om]\cap Q.
  \end{align*}
\end{itemize}
\end{dnt}

Before quoting the solvability results in \cite{JAC17,EF15} for system \eqref{abeqs1}-\eqref{abeqs3}, we need to further introduce the set of \textit{subsolutions}:
\begin{align*}
X_0
=\Big\{\bv \big| \bv\in C_{\textrm{weak}}([0,T];L^2(\Om;\R^d))\cap L^{\infty}((0,T)\times\Om;\R^d),
 \bv(0,\cd)=\bv_0, \bv(T,\cd)=\bv_T,
\end{align*}
\begin{align*}
\p_t\bv  +  \di_x \mathbb{F} = 0,  \di_x\bv = 0\,\, \textrm{in} \,\, \mathcal{D}'((0,T)\times\Om;\R^d), \,\, \textrm{for}\,\, \textrm{some}\,\,
 \bv \in C(Q;\R^d),
\end{align*}
\begin{align*}
\mathbb{F}\in L^{\infty}((0,T)\times\Om;\R_{\textrm{sym},0}^{d\times d})  \cap C(Q;\R_{\textrm{sym},0}^{d\times d})
\end{align*}
\begin{align*}
\sup_{(t,x)\in Q,t>\tau} \f{d}{2}\lambda_{\textrm{max}}\left[\f{(\bv + \bh[\bv])\otimes(\bv + \bh[\bv])}{r[\bv]}  -   \mathbb{F}  +  \mathbb{H}[\bv]\right]
  -   e[\bv]   <0\,\,  \textrm{for} \,\, \textrm{any}\,\, 0<\tau<T \Big\},
\end{align*}
where $\lambda_{\textrm{max}}[A]$ denotes the maximal eigenvalue of a symmetric matrix $A$.
Now, we can state the following existence result for \eqref{abeqs1}-\eqref{abeqs3}, see  \cite{JAC17,EF15}:
\begin{prop}\label{prop2}
Let the operators $\bh$, $r$, $\mathbb{H}$ and $e$ be $Q$-continuous, where $Q\subset [(0,T)\times\Om]$ is an open set satisfying $|Q|=|(0,T)\times\Om|$.
In addition, suppose that $r[\bv]>0$ and that the mapping $\bv\mapsto 1/r[\bv]$ is continuous in the sense specified in Definition \ref{def3}. Finally, assume that the set of subsolutions $X_0$ is non-empty and bounded in $L^{\infty}((0,T)\times\Om;\R^d)$. Then, problem \eqref{abeqs1}-\eqref{abeqs3} admits infinitely many solutions.
\end{prop}


\subsection{Recast \eqref{Eqs-solv}-\eqref{initial-solv} into the abstract Euler form}
In order to apply Proposition \ref{prop2} to prove the solvability of \eqref{Eqs-solv}-\eqref{initial-solv}, we need to firstly recast them into the form of \eqref{abeqs1}-\eqref{abeqs3}. If we can further verify that assumptions in Proposition \ref{prop2} hold, then existence of solutions for the system \eqref{Eqs-solv}-\eqref{initial-solv} is proven. To this end, we take $Q=(0,T)\times\Om$.

\subsubsection{Momentum decomposition and kinetic energy}
Following \cite{JAC17} one can write the momentum $\rho\bu$ in the form
\begin{align*}
\rho\bu = \bv  +  \bV   +   \na_x\Psi,
\end{align*}
where
\begin{align*}
\di_x\bv=0, \quad \int_{\Om}\Psi(t,\cd)dx=0, \quad \int_{\Om}\bv(t,\cd)dx=0, \quad \bV=\bV(t)\in \R^d.
\end{align*}
Similarly, we write the initial momentum $\rho_0\bu_0$ as
\begin{align*}
\rho_0\bu_0=\bv_0 + \bV_0 + \na_x\Psi_0, \quad \di_x\bv_0=0, \quad \int_{\Om}\bv_0dx=\int_{\Om}\Psi_0dx=0, \quad \bV_0=\f{1}{|\Om|}\int_{\Om}\rho_0\bu_0dx.
\end{align*}
Accordingly, we may fix $\rho\in C^2([0,T]\times\Om)$ such that for a certain potential $\Psi$,
\begin{align*}
\p_t\rho  +   \Del_x\Psi  =  0  \quad \textrm{in} \,\, (0,T)\times\Om,
\end{align*}
\begin{align*}
\p_t\rho(0,\cd) = -\Del_x\Psi_0, \quad \Psi(0,\cd) = \Psi_0, \quad \int_{\Om}\Psi(t,\cd)dx = 0 \quad \textrm{for} \,\, \textrm{any} \,\, t\in[0,T].
\end{align*}
Hence, in the sequel, we assume that that
\begin{align*}
\rho\in C^2([0,T]\times\Om), \quad \Psi\in C^1([0,T];C^3(\Om))
\end{align*}
are fixed functions. Based on the above decomposition, equation \eqref{Eqs-solv} reduces to
\begin{align}\label{Eqs-de1}
&\p_t\bv  +  \p_t\bV   +   \di_x\left(\f{(\bv+\bV+\na_x\Psi)\otimes(\bv+\bV+\na_x\Psi)}{\rho}     +    (p(\rho)  +  \p_t\Psi)\mathbb{I}\right)  \nn \\
&\quad \quad =  -(\na_x K\ast\rho)\rho    -  (\bv + \bV + \na_x\Psi)   -   \rho\na_x\Phi,   \nn \\
&\di_x\bv = 0.
\end{align}
In order to match \eqref{abeqs2}, we fix the "kinetic energy" so that
\begin{align}\label{Eqs-de2}
\f12\f{|\bv + \bV + \na_x\Psi|^2}{\rho}  =  e  \equiv \Pi  -  \f{d}{2}(p(\rho) + \p_t\Psi),
\end{align}
where $\Pi=\Pi(t)$ is a spatially homogeneous function to be determined later. Substituting \eqref{Eqs-de2} into \eqref{Eqs-de1}, one can therefore rewrite \eqref{Eqs-de1} as
\begin{align}\label{Eqs-de3}
&\p_t\bv  +  \p_t\bV   +   \di_x\left(\f{(\bv+\bV+\na_x\Psi)\odot(\bv+\bV+\na_x\Psi)}{\rho}  \right)  \nn \\
&\quad \quad =  -(\na_x K\ast\rho)\rho    -  (\bv + \bV + \na_x\Psi)   -   \rho\na_x\Phi  := \mathfrak{E} ,   \nn \\
&\di_x\bv = 0.
\end{align}


\subsubsection{Fix $\bV$ and recast \eqref{Eqs-de3} into abstract form}
One can easily notice from \eqref{Eqs-de3} that there are still two unknowns $\bv$ and $\bV$. So our first goal in this subsubsection is to fix $\bV$ so that
\eqref{Eqs-de3} can be converted to a "balance law" with a source term of zero mean. To this end, solving the following ODE:
\begin{align*}
\f{d\bV}{dt}  +   \bV
=  -  \f{1}{|\Om|}\int_{\Om}(\na_x K\ast\rho)\rho dx   -   \f{1}{|\Om|}\int_{\Om}\rho\na_x\Phi dx
\end{align*}
with initial data $\bV(0)=\bV_0$, one can see that $\bV=\bV[\bv]$ depends linearly on the fixed function $\rho$. Thus, we can therefore rewrite \eqref{Eqs-de3} as
\begin{align}\label{Eqs-de4}
\p_t\bv   +   \di_x\left(\f{(\bv+\bV+\na_x\Psi)\odot(\bv+\bV+\na_x\Psi)}{\rho}  \right)   &=   \mathfrak{E}  -  \f{1}{|\Om|}\int_{\Om}\mathfrak{E}dx,  \nn \\
\di_x\bv &= 0.
\end{align}
Obviously, the expression on the right-hand-side of \eqref{Eqs-de4} has zero integral mean at any time $t$. Hence, referring \cite{JAC17} for more details, we can find a vector $\bw=\bw[\bv]$ satisfying
\begin{align*}
-\di_x\left(\na_x\bw + \na_x^{\top}\bw - \f{2}{d}\di_x\bw\mathbb{I}\right)
= \mathfrak{E}  -  \f{1}{|\Om|}\int_{\Om}\mathfrak{E}dx\quad \textrm{in} \,\, \Om\,\, \textrm{for} \,\, \textrm{any}\,\, \textrm{fixed} \,\, t\in[0,T].
\end{align*}
Denoting
\begin{align}\label{H}
\mathbb{H}[\bv]=\na_x\bw + \na_x^{\top}\bw - \f{2}{d}\di_x\bw\mathbb{I},
\end{align}
one can thus transform system \eqref{Eqs-solv}-\eqref{initial-solv} to the form coincide with \eqref{abeqs1}-\eqref{abeqs3}, namely:

Find a vector field $\bv\in C_{\textrm{weak}}([0,T];L^2(\Om;\R^d))$ satisfying
\begin{align*}
\p_t\bv  +  \di_x\left(\f{(\bv + \bV[\bv] +  \na_x\Psi)\odot(\bv + \bV[\bv] + \na_x\Psi)}{\rho}  +   \mathbb{H}[\bv]\right)   =   0,
\quad \di_x\bv = 0
\end{align*}
in $\mathcal{D}'((0,T)\times\Om;\R^d)$,
\begin{align}\label{abeqs5}
\f12\f{|\bv + \bV[\bv] + \na_x\Psi|^2}{\rho}  =   e[\bv]
\equiv \Pi - \f{d}{2}(p(\rho) + \p_t\Psi) \quad \textrm{for} \,\, \textrm{a.a.}\,\, (t,x)\in (0,T)\times\Om,
\end{align}
\begin{align*}
\bv(0,\cd)=\bv_0, \quad \bv(T,\cd)=\bv_T.
\end{align*}


\subsection{Proof of Theorem \ref{thm2}}
Taking $r[\bv]=\rho$, $\bh[\bv]=\bV[\bv] + \na_x\Psi$, $\mathbb{H}[\bv]$ defined by \eqref{H}, and $e[\bv]$ defined by \eqref{abeqs5}, one can easily verify that they are $Q$-continuous. Then Theorem \ref{thm2} can be proved if we are able to show that $X_0$ is non-empty and bounded in $L^{\infty}((0,T)\times\Om;\R^d)$.

To this end, taking $\bv_T=\bv_0$, $\bv=\bv_0$, and $\mathbb{F}=0$ in the definition of $X_0$ and choosing $\Pi=\Pi(t)$ to be large enough satisfying
\begin{align*}
\sup_{(t,x)\in Q,t>\tau} \f{d}{2}\lambda_{\textrm{max}}\left[\f{(\bv_0 + \bV[\bv_0] + \na_x\Psi)\otimes(\bv_0 + \bV[\bv_0] + \na_x\Psi)}{\rho}   +  \mathbb{H}[\bv_0]\right]
  -   \Pi(t)  +  \f{d}{2}(p(\rho) + \p_t\Psi)<0
\end{align*}
for any $0<\tau<T$, one can claim that there exists $\Pi_0>0$ such that the above inequality holds whenever $\Pi(t)\ge\Pi_0$ for all $t\in[0,T]$.
Consequently, $\bv_0\in X_0$ and thus $X_0$ is non-empty.

In order to prove $X_0$ is bounded in $L^{\infty}((0,T)\times\Om;\R^d)$, we firstly recall the purely algebraic inequality \cite{DL10},
\begin{align}\label{alg-in}
\f12\f{|M|^2}{\tilde r}  \le \f{d}{2}\tilde\lambda_{\textrm{max}}\left[\f{M\otimes M}{\tilde{r}}  -  \mathbb{H}\right] \quad
\textrm{whenever}\,\, \mathbb{H}\in\R_{\textrm{sym},0}^{d\times d},\ M\in \R^d, \tilde{r} \in (0,\infty).
\end{align}
Fixing $\Pi(t)$ according to the above discussions, for any $\bv\in X_0$, we have by using the definition of $X_0$
\begin{align*}
\f{d}{2}\lambda_{\textrm{max}}\left[\f{(\bv + \bV + \na_x\Psi)\otimes(\bv + \bV + \na_x\Psi)}{\rho}   -   (\mathbb{F} - \mathbb{H}[\bv])\right]
<\Pi(t)  -  \f{d}{2}(p(\rho) + \p_t\Psi).
\end{align*}
By the definition of $\mathbb{H}$ in \eqref{H}, one can obtain that $\mathbb{H}[\bv]\in \R_{\textrm{sym},0}^{d\times d}$. Applying the inequality \eqref{alg-in}, we have
\begin{align*}
\f12\f{|\bv + \bV + \na_x\Psi|^2}{\rho}   <  \Pi(t)  -  \f{d}{2}(p(\rho) + \p_t\Psi),
\end{align*}
which implies that $X_0$ is bounded in $L^{\infty}((0,T)\times\Om;\R^d)$. So far, all the assumptions in Proposition \ref{prop2} hold, and the proof of Theorem \ref{thm2} directly follows now by using Proposition \ref{prop2}.


\section*{Acknowledgement}
JAC was partially supported by the EPSRC grant number EP/P031587/1. YP was supported by the CSC (China Scholarship Council) during her visit at Imperial College London. AWK was partially supported by the grant Iuventus Plus no. 0871/IP3/2016/74 of Ministry of Sciences and Higher Education RP. YP's vist at Institute of Mathematics, Polish Academy of Sciences was supported by the programme "Guests of IMPAS".

 {\small

(Jos\'{e} A. Carrillo)\\
DEPARTMENT OF MATHEMATICS, IMPERIAL COLLEGE LONDON, LONDON SW7 2AZ, UNITED KINGDOM

E-mail address: carrillo@imperial.ac.uk
\smallskip\\

(Yingping Peng)\\
SCHOOL OF MATHEMATICAL SCIENCES, UNIVERSITY OF
ELECTRONIC SCIENCE AND TECHNOLOGY OF CHINA, CHENGDU 611731, CHINA \\
AND\\
DEPARTMENT OF MATHEMATICS, IMPERIAL COLLEGE LONDON, LONDON SW7 2AZ, UNITED KINGDOM

E-mail address: yingping$\_$peng@163.com
\smallskip\\

(Aneta Wr\'{o}blewska-Kami\'{n}ska)\\
INSTITUTE OF MATHEMATICS, POLISH ACADEMY OF SCIENCES, \'{S}NIADECKICH 8, 00-656 WARSZAWA, POLAND

E-mail address: awrob@impan.pl


\begin{thebibliography}{99}

\bibitem{PA-2} P. Antonelli and P. Marcati,
 \it The quantum hydrodynamics system in two space dimensions,
 \rm Arch. Rational Mech. Anal, 203(2012), 499-527.


\bibitem{PA-1} P. Antonelli and P. Marcati,
 \it On the finite energy weak solutions to a system in quantum fluid dynamics,
 \rm Comm. Math. Physics, 287(2009), 657-686.

\bibitem{AB} A. Blanchet, J. A. Carrillo and P. Laurencot,
 \it Critical mass for a Patlak-Keller-Segel model with degenerate diffusion in higher dimensions,
 \rm Calc. Var., 35(2009), 133-168.


\bibitem{DB} D. Brandon, T. Lin and R. C. Rogers,
 \it Phase transitions and hysteresis in nonlocal and order-parameter modles,
 \rm Meccanica, 30(1995), 541-565. Microstructure and phase transitions in solods (Udine, 1994).


\bibitem{CH14} E. Chiodaroli,
 \it A counterexample to well-posedness of energy solutions to the compressible Euler system,
 \rm J. Hyperbolic Differ. Eqs., 11(2014), 493-519.


\bibitem{JAC-1} J. A. Carrillo, Y.-P. Choi,
 \it Quantitative error estimates for the large friction limit of Vlasov equation with nonlocal forces,
 \rm arXiv: 1901.07204v1 [math. AP].


\bibitem{VC} V. Calvez, J. A. Carrillo, F. Hoffmann,
 \it Equilibria of homogeneous functionals in the fair-competition regime,
 \rm Nonlinear Anal., 159(2017), 85-128.

\bibitem{CCT} J. A. Carrillo, Y.-P. Choi, O. Tse,
 \it Convergence to equilibrium in Wasserstein distance for damped Euler equations with interaction forces,
  \rm Comm. Math. Phys., 365(2019), no. 1, 329-361.

\bibitem{JAC17} J. A. Carrillo, E. Feireisl, P. Gwiazda and A. \'{S}wierczewska-Gwiazda,
 \it Weak solutions for Euler systems with non-local interactions,
 \rm J. London Math. Soc., 95(2017), no. 2, 705-724.

\bibitem{CH15} E. Chiodaroli, E. Feireisl and O. Kreml,
 \it On the weak solutions to the equations of a compressible heat conducting gas,
 \rm Ann. Inst. H. Poincar\'{e} Anal. Non Lin\'{e}aire, 32(2015), 225-243.

\bibitem{CG} J.-F. Coulombel, T. Goudon,
 \it The strong relaxation limit of the multidimensional isothermal Euler equations,
 \rm Trans. Am. Math. Soc., 359(2007), 637-648.

 \bibitem{JAC} J. A. Carrillo, F. Hoffmann, E. Mainini and B. Volzone,
 \it Ground states in the diffusion-dominated regime,
 \rm Calc. Var., (2018), 57:127.


\bibitem{JW} L. Chen, L. Hong and J. Wang,
 \it Parabolic elliptic type Keller-Segel system on the whole space case,
 \rm Discrete Contin. Dyn. Syst., 36(2016), no. 2, 1061-1084.

\bibitem{CMV2} J. A. Carrillo, R. J. McCann, C. Villani,
 \it Contractions in the 2-Wasserstein length space and thermalization of granular media.
 \rm Arch. Ration. Mech. Anal. 179(2006), 217-263.

\bibitem{CMV1} J. A. Carrillo, R. J. McCann, C. Villani,
 \it Kinetic equilibration rates for granular media and related equations: entropy dissipation and mass transportation estimates,
  \rm Rev. Mat. Iberoamericana, 19(2003), 971-1018.

 \bibitem{CMD-1} C. M. Dafermos,
 \it  Stability of motions of thermoelastic fluids,
 \rm J. Thermal Stresses,  2(1979), 127-134.



\bibitem{CMD-2} C. M. Dafermos,
 \it The second law of thermodynamics and stability,
 \rm Arch. Ration. Mech. Anal., 70(1979), 167-179.

\bibitem{DL10} C. DeLellis and L. Sz\'{e}kelyhidi,
 \it On admissibility criteria for weak solutions of the Euler equations,
 \rm Arch. Ration. Mech. Anal., 195(2010), 225-260.

 \bibitem{RDJ} R. J. Diperna,
 \it Uniqueness of silutions to hyperbolic conservation laws,
 \rm Indiana Univ. Math. J., 28(1979), 137-188.


\bibitem{DD15} D. Donatelli, E. Feireisl and P. Marcati,
\it Well/ill posedness for the Euler-Korteweg-Poisson system and relatedd problems,
\rm Commun. Partial Differential Equations, 40(2015), 1314-1335.


\bibitem{JED} J. E. Dunn and J. Serrin,
 \it On the thermomechanics of interstitial working,
 \rm Arch. Rational Mech. Anal, 88(1985), 95-133.


\bibitem{EF15} E. Feireisl,
 \it Weak solutions to problems involving inviscid fluids,
 \rm In: Shibata Y., Suzuki Y. (eds) Mathematical Fluid Dynamics, Present and Future. Springer Proceedings in Mathematics \& Statistics, vol 183. Springer, Tokyo. (2016),  377-399.


\bibitem{EF16} E. Feireisl, P. Gwiazda and A. \'{S}wierczewska-Gwiazda,
 \it On weak solutions to the 2d Savage-Hutter model of the motion of a gravity driven avalanche flow,
 \rm Comm. Partial Differential Equations, 41(2016), no. 5, 759-773.


\bibitem{JG} J. Giesselmann, C. Lattanzio and A. E. Tzavaras,
 \it Relative energy for the Korteweg theory and related Hamiltonian flows in gas dynamics,
 \rm Arch. Rational Mech. Anal, 223(2017), no. 3, 1427-1484.

\bibitem{GPK} B. D. Goddard, G. A. Pavliotis, S. Kalliadasis,
 \it The overdamped limit of dynamic density functional theory: rigorous results,
 \rm Multiscale Model. Simul. 10(2012), 633-663.

\bibitem{H1} F. Huang, P. Marcati, R. Pan,
 \it Convergence to the Barenblatt solution for the compressible Euler equations with damping and vacuum,
 \rm Arch. Ration. Mech. Anal., 176(2005), 1-24.

\bibitem{H2} F. Huang, R. Pan,
 \it Asymptotic behavior of the solutions to the damped compressible Euler equations with vacuum,
 \rm J. Differ. Eqs., 220(2006), 207-233.

\bibitem{H3} F. Huang, R. Pan, Z. Wang,
 \it $L^1$ convergence to the Barenblatt solution for compressible Euler equations with damping,
 \rm Arch. Ration. Mech. Anal., 200(2011), 665-689.

\bibitem{QJ} S. Jiang, Q. Ju, H. Li and Y. Li,
 \it Quasi-neutral limit of the two-fluid Euler-Poisson system,
 \rm Commum. Pure Appl. Anal., 9(2010), 1577-1590.

\bibitem{JR} S. Junca, M. Rascle,
 \it Strong relaxation of the isothermal Euler system to the heat equation,
 \rm Z. Angew. Math. Phys., 53(2002), 239-264.

\bibitem{TL} T. Luo and J. Smoller,
 \it Existence and non-linear stability of rotating star solutions of the compressible Euler-Poisson equations,
 \rm Arch. Ration. Mech. Anal., 191(2009), 447-496.


 \bibitem{CL-2} C. Lattanzio, A. E. Tzavaras,
 \it From gas dynamics with large friction to gradient flows describing diffusion theories,
 \rm Comm. Partial Differential Equations, 42(2017), no. 2, 261-290.


 \bibitem{CL-1} C. Lattanzio, A. E. Tzavaras,
 \it Relative entropy in diffusive relaxation,
 \rm SIAM J. Math. Anal., 45(2013), 1563-1584.

 \bibitem{MM} P. Marcati, A. Milani,
 \it The one-dimensional Darcy's law as the limit of a compressible Euler flow,
 \rm J. Differ. Eqs., 84(1990), 129-147.

\bibitem{YP} Y. Peng and Y. Wang,
 \it Convergence of compressible Euler-Poisson equations to incompressible type Euler equations,
 \rm Asymptotic Anal., 41(2005), 141-160.

\bibitem{XR} X. Ren and L. Truskinovsky,
 \it Finite scale microstructures in nonlocal elasticity,
 \rm J. Elasticity, 59(2000), 319-355.

\bibitem{YS} Y. Sugiyama,
 \it Global existence in sub-critical cases and finite time blow-up in super-critical cases to degenrate Keller-Segel systems,
 \rm Differential Integral Equations, 9(2006), 841-876.

\bibitem{V} C. Villani,
 \it Topics in optimal transportation,
 \rm Graduate Studies in Mathematics, 58. American Mathematical Society, Providence, RI, 2003. xvi+370 pp.


\end{thebibliography}
\end{document}